\documentclass{article}
\title{A poor man's improvement on Zhang's result: there are infinitely many prime gaps less than 60 million}
\author{T. S. Trudgian\footnote{Supported by ARC Grant DE120100173.} \\
Mathematical Sciences Institute\\
The Australian National University\\
timothy.trudgian@anu.edu.au}
\usepackage{comment}
\usepackage{url}
\usepackage{amsthm}
\usepackage{amsmath}
\usepackage{amssymb}
\begin{document}
\maketitle

%\begin{abstract}
%In this article the author takes the easiest possible path in improving Zhang's result.
%\end{abstract}
\noindent
Consider a set $\mathcal{H} = \{ h_{1}, h_{2}, \ldots, h_{k_{0}}\}$, composed of distinct, non-negative integers, in which $k_{0}$ denotes a natural number. Call the set $\mathcal{H}$ \textit{admissible} if $\nu_{p}(\mathcal{H})<p$ for every prime $p$, where $\nu_{p}(\mathcal{H})$ is the number of distinct residue classes modulo $p$ that are covered by the elements $h_{i}$.

In \cite[Thm.\ 1]{Zhang} Zhang shows that if $k_{0}= 3.5\times 10^{6}$ and $\mathcal{H}$ is admissible, then there are infinitely many positive integers $n$ for which $\{n+ h_{1}, \ldots, n+ h_{k_{0}}\}$ contains at least two primes. He follows this with three lines of deduction to prove the remarkable result

$$\liminf_{n\rightarrow \infty} (p_{n+1} - p_{n}) < 7 \times 10^{7}.$$
It is with this three-line deduction --- the infinitely easier portion of Zhang's paper --- that this article is concerned. Given a value of $k_{0}$, how can one find an admissible set $\mathcal{H}$ for which the length of $\mathcal{H}$, defined to be $h_{k_{0}} - h_{1}$, is small? 

\subsubsection*{An upper bound for the length of $\mathcal{H}$}
Consider a set of $k_{0}$ primes 
$$\mathcal{H} = \{p_{m+1}, p_{m+2}, \ldots, p_{m+ k_{0}}\},$$
where $m$ is a non-negative integer to be determined momentarily; Zhang takes $m= k_{0}$. When is $\mathcal{H}$ an admissible set?

Consider primes $p\leq p_{m}$. Since $p_{m+i} \not\equiv\, 0 \,(\mbox{mod } p) $ for all $1\leq i \leq k_{0}$, it follows that $\nu_{p}(\mathcal{H})\leq p-1< p$ for all $p\leq p_{m}$.

Now consider primes $p\geq p_{m+1}$. We should like to show that there is not enough room in the set to fill all of the residue classes modulo $p$. Therefore for primes $p\in \mathcal{H}$ we wish to show that $k_{0} - 1 < p -1$. It is sufficient to show that $p_{m+1} > k_{0}.$ A quick computation shows that one may choose $m= 250,150$. 
Therefore the maximal gap between primes in $\{n + p_{m+1}, \ldots, n + p_{m+ k_{0}}\}$ is $p_{m+ k_{0}} - p_{m+1} = 59,874,594.$

In 2013, between 30th May and 3rd June, a considerable amount of work was undertaken by Morrison, Tao, et al. \cite{SBS} which not only improved on the method of exhibiting small gaps, but also improved on the value of $k_{0}$. To date, the smallest permissible value of $k_{0}$ is $341,640$, which leads to a prime gap not exceeding $4,802,222$.

\subsubsection*{An lower bound for the length of $\mathcal{H}$}

The set  $\{1, 2\ldots, k_{0}\}$ cannot be admissible since, \textit{inter alia}, both even and odd numbers are present. Therefore, at the very least, we must impose that our set be of the form $\{r_{1}, r_{1} + 2, \ldots, r_{1} + 2(k_{0} - 1)\}$, where $r_{1}$ is either 1 or 0  modulo 2. Such a set has $k_{0}$ elements, and length $2(k_{0} - 1)$. It may be that this set is not admissible; the point to note is that the \textit{minimal} length of an admissible set must be bounded below by $2(k_{0} - 1)$.

We may generalise this approach by noting that we can fill, at the most, $p_{i}-1$ residue classes modulo $p_{i}$, for each $p_{i}$. At best, we may include $R_{m} = (p_{1} - 1)\cdots (p_{m} - 1)$ integers modulo $M_{m} =p_{1}\cdots p_{m}$, where $m$ is an integer that we shall determine momentarily. For $1\leq i \leq j \leq R_{m}$ let $r_{i}\leq r_{j}$ run through the $R_{m}$ residues modulo $M_{m}$. Consider the sets
$$T=\{r_{1}, r_{2}, \ldots, r_{R_{m}}, \ldots, r_{1} + (a-1)M_{m}, r_{2} + (a-1)M_{m}, \ldots, r_{R_{m}} + (a-1)M_{m}\}$$
and
$$ T' = \{ r_{1} + a M_{m}, \ldots, r_{n} + a M_{m}\},$$
in which $1\leq n \leq R_{m}$ and in which $a$ is chosen such that $|T|< k_{0}$ and $|T \cup T'| \geq k_{0}$.
It follows that $R_{m}a < k_{0} \leq R_{m}(a+1)$, whence the length of $T \cup T'$ bounded below by
\begin{equation}\label{klo}
a M_{m} \geq M_{m} \left(\frac{k_{0}}{R_{m}} - 1\right).
\end{equation}
Given a value of $k_{0}$ one may choose the value of $m$ maximising the right-side of (\ref{klo}). When $k_{0} = 341,640$, one should choose $m=6$, which gives a gap at least as large as 1,751,112.
\begin{comment}

\subsubsection*{Method B}

The set  $\{1, 2\ldots, k_{0}\}$ cannot be admissible since, \textit{inter alia}, both even and odd numbers are present. Delete elements congruent to $0 \,(\mbox{mod } 2)$, and expand the set to form the new set $\{1, 3, \ldots, 2k_{0}-1\}$ which has $k_{0}$ integers. By deleting elements congruent to $0 \,(\mbox{mod } 3)$, the set consists of elements congruent to $1, 5\,(\mbox{mod } 6)$.  We continue, as required, this approach of expansion and deletion of members divisible by primes.

Consider $T_{0} = k_{0}$ and $T_{m} = T_{m-1} \frac{p_{m}}{p_{m} - 1}$. Now define 
$$S_{m} = \{n:\,  n\leq \lfloor T_{m}\rfloor \; \textrm{and}\; (n, p_{1}\cdots p_{m}) = 1\},$$
so that $S_{m}$ consists of at least $k_{0}$ integers that are not divisible by any of $p_{1}, \ldots, p_{m}$. Once $p_{m+1}^{2}> T_{m}$, or, equivalently,
\begin{equation}\label{impo}
k_{0} \prod_{r=1}^{m} \frac{p_{r}}{p_{r} - 1} < p_{m+1}^{2},
\end{equation}
we can stop the process. Whereas none of the primes $p_{1}, \ldots, p_{m}$ divide elements in $S_{m}$, (\ref{impo}) shows that, for all $i\geq 1$ no prime $p_{m+i}$ divides elements in $S_{m}$.

Therefore, for the least value of $m$ that satisfies (\ref{impo})
we have an admissible set $\{1, \ldots, T_{m}\}$ and hence a maximal gap of $T_{m} -1.$
Solving (\ref{impo}) for $k_{0} = 3.5\times 10^{6}$ gives $m=944$ and a maximal gap of 55,622,435.
\end{comment}
\subsubsection*{Comparison of bounds}

These bounds appear to be wasteful. Consider, for example, the data in \cite[p.\ 832]{GPY}. There, conditional values of $k_{0}$ are given along with the corresponding minimal length of the $k_{0}$-tuple. Table 1 compares the results in \cite{GPY}, obtained by an exhaustive computational search, with the upper and lower bounds obtained here.

\begin{table}[ht]
\caption{Comparison of bounds for gaps between successive primes}
\label{optable}
\centering
\begin{tabular}{c c c c}
\hline\hline
 $k_{0}$ & Upper bound & Lower bound & Length in \cite[p.\ 829, 832]{GPY}  \\[2 ex] 
6 & 16 & 12 &16\\
10 & 32 & 24 & 32\\
12  & 46 & 30 & 42\\
%32 & 162 & 90& 158\\
65  & 364 & 189 & 336\\
%111 & 696 &387& 634\\
193  & 1292 &694& 1204\\
1000 & 8424 & 4165 & \\
10000 & 109152 & 45815 &\\
$341,640$ & $5,005,362^{*}$ & $1,751,112$ & \\
 \hline\hline
\end{tabular}

%These points roughly fit a straight line of equation $y = 0.925x + 3.04$, where $y$ denotes the fourth column above, and $x$ denotes the third column. Throwing all caution into the wind one could speculate that the bound of $59,874,594$ therefore corresponds to an `actual' gap of the order of $5.5\times 10^{7}$.
\begin{footnotesize}*Note that, as mentioned on page 2, this has been improved to 4,802,222.
\end{footnotesize}
\end{table}
While the upper bound gives the correct answer for small values of $k_{0}$, it  becomes increasingly profligate as $k_{0}$ increases; the lower bound appears to be ubiquitously impotent.

The method of searching by brute force, potentially another poor man's improvement, appears to be next to hopeless. For, given $k_{0}$ distinct non-negative integers of size at most $N$, there are $\binom{N}{k_{0}}$ possible $k_{0}$-tuples. Even with the modest value of $N= 7\times 10^{6}$ one faces the daunting prospect of searching for an admissible ($341,640$)-tuple amongst more than $10^{2\times10^{5}}$ possible candidates.

%It is with great interest that the author awaits further discussion and research on Zhang's result.

I am grateful for Scott Morrison's providing me with data for $k_{0} = 341,640$, and for the interesting work undertaken by him, Terry Tao, and others in \cite{SBS}.

\bibliographystyle{plain}
\bibliography{themastercanada}

\end{document}